\documentclass[10pt]{article}
\usepackage{amsfonts}
\usepackage{amsmath}
\usepackage{a4}
\usepackage{amssymb}
\usepackage{theorem}
\usepackage{psfrag}

\newtheorem{theorem}{Theorem}

{\theorembodyfont{\rmfamily}}
{\theorembodyfont{\rmfamily}}
\newtheorem{proposition}[theorem]{Proposition}
\numberwithin{equation}{section}
\numberwithin{figure}{section}

\newcommand{\RR}{\mathbb{R}}
\newcommand{\PP}{\mathbb{P}}
\newcommand{\ZZ}{\mathbb{Z}}

\newcommand{\EE}{\mathbb{E}\,}

\newcommand{\cL}{\mathcal{L}}

\newcommand{\cU}{\mathcal{U}}

\newcommand{\Var}{\operatorname{Var}}

\renewcommand{\Box}{\square}

\newcommand{\bu}{\mathbf{u}}

\allowdisplaybreaks

\begin{document}

\title{A universality property for last-passage 
percolation paths close to the axis}
\author{Thierry Bodineau\footnote{Laboratoire de Probabilit\'es et Mod\`eles Al\'eatoires,
Universit\'e Pierre et Marie Curie - Bo{\^\i}te courrier 188
75252 Paris Cedex 05, France.
Email: \texttt{bodineau@gauss.math.jussieu.fr}} \, and 
James B.\ Martin\footnote{LIAFA, Universit\'e Paris 7,
case 7014, 2 place Jussieu, 75251 Paris Cedex 05, France. Email: 
\texttt{martin@liafa.jussieu.fr}}}
\date{}


\maketitle


\begin{abstract}
We consider a last-passage directed percolation model
in $\ZZ_+^2$, with i.i.d.\ weights whose common 
distribution has a finite $(2+p)$th moment.
We study the fluctuations of the passage time from
the origin to the point $\big(n,n^{\lfloor a \rfloor}\big)$.
We show that, 
for suitable $a$ (depending on $p$),
this quantity, appropriately scaled,
converges in distribution as $n\to\infty$
to the Tracy-Widom distribution,
irrespective of the underlying weight distribution.
The argument uses a coupling to a Brownian directed
percolation problem and the strong approximation 
of Koml\'os, Major and Tusn\'ady.
\end{abstract}

\section{Introduction}

The concept of \textit{universality class} plays a key role in statistical
mechanics, making it possible to classify a huge variety of models and
phenomena by means of well chosen scaling exponents.
For instance, many growth models are expected to 
share similar properties which
fall into the framework of the KPZ universality class -- 
see for example the survey by Krug and Spohn \cite{KrugSpohn}.

In this note, we focus on the particular example of 
\textit{directed last-passage percolation}.
Let $\omega_i^{(r)}$, $i\geq 0, r\geq 1$ be i.i.d.\ random variables.
We consider directed paths in the lattice $\ZZ_+^2$,
each step of which increases one of the coordinates by 1.
For $n\geq 0$, $k\geq 1$, 
the \textit{(last-)passage time} to the point $(n,k)$ is
defined by 
\begin{equation}\label{tdef}
T(n,k)
=
\max_{\pi\in\Pi(n,k)}
\left\{\sum_{(i,r)\in\pi} \omega_i^{(r)}\right\},
\end{equation}
where $\Pi(n,k)$ is the set
of directed paths from $(0,1)$ to $(n,k)$.
More precisely,
\begin{multline*}
\Pi(n,k)=
\Bigg\{
\left(z_1,z_2,\dots,z_{n+k}\right)\in\left(\ZZ^2_+\right)^{n+k}:
z_1=(0,1), \;\;
z_{n+k}=(n,k),
\\ 
z_{j+1}-z_j \in \{(0,1),(1,0)\} 
\text{ for } 1\leq j\leq n+k-1
\Bigg\}.
\end{multline*}

When the underlying weight distribution is exponential or geometric,
scaling exponents for this model are rigorously known. 
The deviation from a straight line of the optimal path to the point $(n,n)$ 
is of the order $n^{2/3}$ (corresponding to the
exponent $\xi=2/3$) \cite{BDMMZ,Johanssontransversal}, while
the fluctuations of the passage time $T(n,n)$ 
are of the order $n^{1/3}$ (corresponding to the exponent $\chi=1/3$).
In fact, one can give much more precise information: for the exponential 
distribution with mean 1, say, it is shown in \cite{Johshape} that
\begin{equation}
\label{TWlimit}
n^{-1/3}\big[T(n,n)-4n\big]\to F_{\textnormal{TW}},
\end{equation}
where $F_{\textnormal{TW}}$ is the ``Tracy-Widom'' distribution, 
which also appears as the asymptotic distribution of the largest eigenvalue 
of a GUE random matrix.
It is expected (but not yet proved) 
that the same scaling exponents, and indeed
the same asymptotic distribution in (\ref{TWlimit}), should
hold for a wide class of underlying weight distributions.

In this note we give a universality result for the quantities
$T\left(n,\lfloor n^a \rfloor\right)$ for $a<1$.
Thus, we are concerned with passage times to points 
which are asymptotically rather close to the horizontal axis,
for a general class of underlying weight distribution.
Our result is the following:
\begin{theorem}
\label{convergencetheorem}
Suppose that $\EE |\omega_i^{(r)}|^p$ for some $p>2$,
with $\mu=\EE \omega_i^{(r)}$, $\sigma^2=\Var(\omega_i^{(r)})$.
Then for all $a<\frac{6}{7}\left(\frac{1}{2}-\frac{1}{p}\right)$,
\begin{equation}
\label{convergence}
\frac{T\left(n,\lfloor n^a \rfloor \right)-n\mu-2 \sigma \, n^{\frac{1+a}{2}}}
{\sigma n^{\frac{1}{2}-\frac{a}{6}}}
\to F_{\textnormal{TW}}
\end{equation}
in distribution. 
In particular, if the weight distribution has finite moments of all orders,
then (\ref{convergence}) holds for all $a<3/7$.
\end{theorem}

Heuristically the theorem can be understood as follows.
As the optimal path goes from the origin to $(n,n^{\lfloor a\rfloor})$, 
one can imagine that between each step upwards,
the path typically takes on the order of $n^{1-a}$ steps to the right.
Thus it should behave as the optimal path from the origin to $(n^a,n^a)$
in a last percolation model with Gaussian weights of variance $n^{1-a}$.
On the renormalized scale the expected fluctuations are of order
$(n^a)^{1/3}$.
In this way, we recover the fluctuation exponent
\[
\hat \chi= \frac{1-a}{2} + \frac{a}{3} = \frac{1}{2}-\frac{a}{6} \, .
\]

This heuristic is made precise by coupling the discrete model 
with a Brownian directed percolation model 
for which the fluctuations have been
explicitly computed. This is done using the strong approximation 
of a random walk by a Brownian motion due to Koml\'os, Major and Tusn\'ady.
(We note that this strong approximation has already been applied
to similar last-passage percolation models by Glynn and Whitt
\cite{GlyWhi}).

A different sort of universality result for paths near the axis 
in directed percolation models is given in \cite{jbmshape}.
By subadditivity, one has the 
convergence $n^{-1}T(n,\lfloor xn \rfloor)\to\gamma(x)$ 
a.s.\ and in $\cL_1$, for some function $\gamma$.
Under the hypothesis of Theorem \ref{convergencetheorem},
it's shown that
$\gamma(x)=\mu+2\sigma\sqrt{x} + o(\sqrt{x})$ as $x\downarrow0$.

The proof of Theorem \ref{convergencetheorem} is given in the next
section. In Section \ref{further},
we make some comments on related models and possible extensions.

\section{Fluctuations of the passage-time}

We first introduce
the \emph{Brownian directed percolation model}.
Let $B^{(r)}_t$, $t\geq 0$, $r\geq 1$ be a sequence of
independent standard Brownian motions.
For $t>0$, $k\geq 1$, define
\[
\cU(t,k)=
\left\{(u_0,u_1,\dots,u_k)\in\RR^{k+1}:
0=u_0\leq  u_1\leq\dots\leq u_k=t
\right\},
\]
and then let
\begin{equation}
\label{ldef}
L(t,k)=\sup_{\bu\in\cU(t,k)}\sum_{r=1}^k
\left[
B^{(r)}_{u_r}-B^{(r)}_{u_{r-1}}
\right].
\end{equation}

One can rewrite the definition of $T(n,k)$ at (\ref{tdef}) in
an analogous way:
\begin{equation}
\label{tdef2}
T(n,k)=\sup_{\bu\in\cU(n,k)}\sum_{r=1}^k
\left[
S^{(r)}_{\lfloor u_r\rfloor+1}-S^{(r)}_{\lfloor u_{r-1}\rfloor}
\right],
\end{equation}
where $S^{(r)}_m=\sum_{i=0}^{m-1} \omega^{(r)}_i$.

The random variable $L(1,k)$ has the same distribution as 
the largest eigenvalue of a $k\times k$ GUE random matrix 
\cite{GUEs}, \cite{GraTraWid}, \cite{OcoYor2}.
Hence in particular (e.g.\ \cite{TracyWidomsurvey})
\[
k^{1/6}\big[L(1,k)-2\sqrt{k}\big]
\to 
F_{\textnormal{TW}}
\]
in distribution, where $F_{\textnormal{TW}}$ is the 
Tracy-Widom distribution.

By Brownian scaling, $L(t,k)$ has the same distribution as $\sqrt{t}L(1,k)$.
Using this we get, for any $0<a\leq 1$,
\begin{equation}
\label{TWconv}
\frac{L\left(n,\lfloor n^a \rfloor \right)-2n^{\frac{1+a}{2}}}
{n^{\frac{1}{2}-\frac{a}{6}}}
\to F_{\textnormal{TW}}
\end{equation}
in distribution.

Theorem \ref{convergencetheorem} says
that the same distributional limit 
as in (\ref{TWconv})
(in particular,
with the same order of fluctuations) 
occurs for the law of $T(n, \lfloor n^a\rfloor)$, for a general underlying
distribution of the weights $\omega^{(r)}_i$, 
if $a$ is sufficiently small.
We will use the following strong approximation 
result, which combines Theorem 2 of 
Major \cite{Major} and Theorem 4 of Koml\'os, Major and Tusn\'ady \cite{KMT2}:
\begin{proposition}
\label{KMTproposition}
Suppose $\omega_i, i=1,2,\ldots$ are i.i.d.\ with
$\EE|\omega_i|^p < \infty$ for some $p>2$, 
and with $\EE \omega_i=0$, $\Var(\omega_i)=1$.
Let $S_m=\sum_{i=0}^{m-1} \omega_i$, $m\geq 1$.

Then there is a constant $C$ 
such that for all $n>0$, there is a coupling of
the distribution of $(\omega_1, \ldots, \omega_n)$ 
and a standard Brownian motion $B_t$, $0\leq t\leq n+1$ 
such that, for all $x\in[n^{1/p}, n^{1/2}]$,
\begin{equation}
\label{strong}
\PP\left(
\max_{m=1,2,\ldots,n+1} \left| B_m-S_m \right|
>x
\right)
\leq Cnx^{-p}.
\end{equation}
\end{proposition}

\noindent\textit{Proof of Theorem \ref{convergencetheorem}}:

We may assume that $\mu=0$ and $\sigma^2=1$,
so that we need to prove that  
\begin{equation}
\label{convergence2}
\frac{T\left(n,\lfloor n^a \rfloor \right)-2n^{\frac{1+a}{2}}}
{n^{\frac{1}{2}-\frac{a}{6}}}
\to F_{\textnormal{TW}}
\end{equation}
in distribution
(for general $\mu$ and $\sigma^2$, 
one can obtain (\ref{convergence}) from (\ref{convergence2})
after replacing $\omega$ by 
$(\omega-\mu)/\sigma$).

If $(\omega^{(r)}_i)_{i,r}$ and $\big(B^{(r)}_t\big)_{t,r}$ are 
all defined on
the same probability space, then from (\ref{ldef}) and (\ref{tdef2})
we get
\begin{multline}
\left|T(n,\lfloor n^a \rfloor)-L(n,\lfloor n^a \rfloor)\right|
\\
\begin{split}
&=
\left|
\sup_{\bu\in\cU(n,\lfloor n^a \rfloor)}\sum_{r=1}^{\lfloor n^a \rfloor}
\left(
S^{(r)}_{\lfloor u_r\rfloor+1}-S^{(r)}_{\lfloor u_{r-1}\rfloor}
\right)
-
\sup_{\bu^\prime \in\cU(n,\lfloor n^a \rfloor)}\sum_{r=1}^{\lfloor n^a \rfloor}
\left(
B^{(r)}_{u^\prime_r}-B^{(r)}_{u^\prime_{r-1}}
\right)
\right|
\\
&\leq
\sup_{\bu\in\cU(n,\lfloor n^a \rfloor)}
\Bigg\{
\sum_{r=1}^{\lfloor n^a \rfloor}
\left|
S^{(r)}_{\lfloor u_r\rfloor +1}
-B^{(r)}_{\lfloor u_r\rfloor +1}
\right|
+
\left|
S^{(r)}_{\lfloor u_{r-1}\rfloor}
-B^{(r)}_{\lfloor u_{r-1}\rfloor}
\right|
\\
&\;\;\;\;\;\;\;\;\;\;\;\;\;\;\;\;\;\;\;\;\;\;\;\;\;\;\;\;\;\;\;
+
\left|
B^{(r)}_{\lfloor u_r\rfloor +1}-B^{(r)}_{u_{r}}
\right|
+
\left|
B^{(r)}_{\lfloor u_{r-1}\rfloor}-B^{(r)}_{u_{r-1}}
\right|
\Bigg\}
\\
&\leq 2\sum_{r=1}^{\lfloor n^a \rfloor}
\left\{
\max_{i=1,2,\ldots, n+1}
\left|
S^{(r)}_i-B^{(r)}_i
\right|
+
\sup_{\substack{0\leq s,t\leq n+1\\|s-t|<2}}
\left|
B^{(r)}_{s}-B^{(r)}_t
\right|
\right\}
\\
\label{VWdef}
&=2\sum_{r=1}^{\lfloor n^a \rfloor}
\left\{
V_n^{(r)} + W_n^{(r)}
\right\},
\end{split}
\end{multline}
where we have defined 
\[
V_n^{(r)}=
\max_{i=1,2,\ldots, n+1}
\left|
S^{(r)}_i-B^{(r)}_i
\right|
\text{ and }
W_n^{(r)}=
\sup_{\substack{0\leq s,t\leq n+1\\|s-t|<2}}
\left|
B^{(r)}_{s}-B^{(r)}_t
\right|.
\]

For each $r=1,\dots,\lfloor n^a \rfloor$
we will couple 
$\left(\omega^{(r)}_0,\omega^{(r)}_1, \ldots, \omega^{(r)}_n\right)$ 
and $B^{(r)}_t$, $0\leq t\leq n+1$ 
as in Proposition \ref{KMTproposition}, maintaining the independence
for different $r$, so that the $V^{(r)}_n$, 
$1\leq r\leq \lfloor n^a \rfloor$ are i.i.d.\ 
with 
\begin{equation}
\label{Vbound}
\PP\left(V^{(r)}_n > x\right)
\leq Cnx^{-p}
\end{equation}
for all $x\in\left[n^{1/p}, n^{1/2}\right]$.

Let $A_1$ be the event 
$\left\{\max_{1\leq r\leq \lfloor n^a\rfloor}
V^{(r)}_n > n^{1/2}
\right\}$.
Then from (\ref{Vbound}),
\begin{align}
\PP(A_1)\leq n^a C n(n^{1/2})^{-p}
= Cn^{a+1-p/2} 
\label{A1to0}
\to 0 \text{ as } n\to\infty,
\end{align}
since by assumption 
$a<\frac67\left(\frac12-\frac1p\right)$
$<p\left(\frac12-\frac1p\right)$
$=\frac{p}2-1$.

Also let $A_2$ be the event 
$\left\{\max_{1\leq r\leq \lfloor n^a\rfloor}
W^{(r)}_n > n^{1/p}
\right\}$.
Using the reflection principle and standard estimates
on the normal distribution,
\begin{align}
\nonumber
\PP(A_2)
&\leq 
n^a \PP\left(
\sup_{\substack{0\leq s,t\leq n+1\\|s-t|<2}}
\left|B^{(1)}_s-B^{(1)}_t\right| > n^{1/p}
\right)
\\
\nonumber
&\leq
n^a
\sum_{i=0}^{n-2}\PP\left(
\sup_{i\leq t\leq i+3} B_t 
-
\inf_{i\leq t\leq i+3} B_t 
> n^{1/p}
\right)
\\
\nonumber
&\leq
n^{a+1}
\PP\left(
\sup_{0\leq t\leq 3}\left|B_t\right|>n^{1/p}/2\right)
\\
\nonumber
&=
4 n^{a+1}\PP\left(B_3>n^{1/p}/2\right)
\\
\nonumber
&\leq c_1 n^{a+1}\exp\left(-c_2 n^{2/p}\right)
\\
\label{A2to0}
&\to 0 \text{ as } n\to\infty.
\end{align}

From (\ref{VWdef}),
(\ref{Vbound}) and the definitions of the events $A_1$ and $A_2$,
we have 
\begin{multline}
\EE\Big[
\left|T(n,\lfloor n^a \rfloor)-L(n,\lfloor n^a \rfloor)\right|
;A_1^C\cup A_2^C
\Big]
\\
\begin{split}
\label{referhere}
&\leq 2n^a \EE\left(V_n^{(1)}+W_n^{(1)};A_1^C\cup A_2^C\right)
\\
&\leq
2n^a\left[
n^{1/p}+\EE\left(
V^{(1)}_n-n^{1/p}; n^{1/p}\leq V^{(1)}_n\leq n^{1/2}\right)+n^{1/p}
\right]
\\
&\leq
2n^a\left(
2n^{1/p}+\int_{n^{1/p}}^{n^{1/2}}\PP\left(V^{(1)}_n>x\right)dx
\right)
\\
&\leq
2n^a\left(
n^{1/p}+\int_{n^{1/p}}^{n^{1/2}}Cnx^{-p}dx
\right)
\\
&=
2n^a\left(
n^{1/p}+C_2n[-x^{-p+1}]^{n^{1/2}}_{n^{1/p}}
\right)
\\
&\leq 
C_3 n^a n^{1/p},
\end{split}\end{multline}
where $C_2$ and $C_3$ are constants independent of $n$.

Together with (\ref{A1to0}) and (\ref{A2to0}), 
this gives that, for any $\epsilon>0$, 
\begin{equation}
\label{Mto0}
\PP\left(\left|T(n,\lfloor n^a \rfloor)-L(n,\lfloor n^a \rfloor)\right|
>n^{a+1/p+\epsilon}\right)
\to 0 \text{ as } n\to\infty.
\end{equation}
The assumption $a<\frac67\left(\frac12-\frac1p\right)$
implies that, for $\epsilon$ sufficiently small, 
$a+\frac1p+\epsilon<\frac12-\frac{a}6$.
Thus
\begin{equation}
\frac{\left|T(n,\lfloor n^a \rfloor)-L(n,\lfloor n^a \rfloor)\right|}
{n^{\frac12-\frac{a}6}}\to 0 \text{ in distribution, as }n\to\infty.
\end{equation}
Using (\ref{TWconv}) we obtain (\ref{convergence2})
as desired.$\hfill\Box$

\medskip

\section{Further remarks}\label{further}

\subsection{Larger values of $a$}
It seems unlikely that the value $3/7$ in Theorem 
\ref{convergencetheorem} represents a real threshold. 
For $a>3/7$, consider the typical difference between
the weight of 
the maximal Brownian path and the weight of
the discrete approximation; this will be large compared to the 
order of fluctuations of the maximum.
However, the standard deviation of this difference
may be smaller; one might expect it to be of order
$n^{a/2}$ rather than order $n^a$, 
since it is composed of $n^a$ terms of constant order
which one expects to become independent as $n$ becomes large.
An argument along these lines would effectively allow us to 
replace $n^a$ by $n^{a/2}$ in (\ref{referhere})
leading to a bound $a<3/4$ rather than $a<3/7$.
However, above $a=3/4$ it seems that the behaviour 
is genuinely different and 
more sophisticated arguments would be required:
the fluctuations of the error in the discrete approximation
to the Brownian path
are likely to be larger than the fluctuations of the maximal
weight itself, and so
one might no longer expect the maximal discrete path
to follow closely the maximal Brownian path (even when 
the weights are ``strongly coupled'' to the Brownian motions as above).

\subsection{Transverse fluctuations}
The exponent $\hat \chi =\frac{1}{2} -
\frac{a}{6}$ which we obtain 
should be related to the transversal fluctuations
of the optimal path away from the straight line $\{ y=n^{a-1} x \}$.
Let us introduce the exponent
\begin{eqnarray*}
\hat \xi = \lim_{n \to \infty}
\frac{1}{2 \log n} \, 
\log \EE \left( \left(v_{\frac{n}{2}} - \frac{n^a}{2}
\right)^2 \right) \, .  \\
\end{eqnarray*}
where, say, $v_i$ is the smallest value $r$ such that
the point $(i,r)$ is contained in the optimal path.

Corresponding to the universal value of $\hat \chi$ obtained,
one would expect that $\hat \xi$
should be also a universal exponent equal to $2a/3$.
To see this, we follow the heuristics explained in the introduction.
On a renormalized level, the optimal path 
should behave as the optimal path from the origin to $(n^a,n^a)$
in a last percolation model with Gaussian weights. 
This would imply that the transverse fluctuations should scale like 
$(n^a)^{2/3}$, where $\xi = 2/3$ is the (predicted) 
standard fluctuation exponent for directed last-passage percolation.

The strategy used in \cite{Johanssontransversal} for the derivation 
of the transverse fluctuations requires not only the
knowledge of the last passage time fluctuations, but also a precise
control of the moderate deviations.
Our approach does not allow us to derive such sharp estimates
(in particular we are missing some uniformity with respect to
the direction of the path). For this reason, 
we do not yet have a proof of the universality of the
transversal fluctuation exponent in our framework.

\subsection{Related models}

The last-passage percolation processes have a natural
interpretation in terms of systems of queues in tandem
(see for example \cite{BBM, GlyWhi}). Considering paths
near the axis corresponds to considering
regimes of very high or very low load in the queueing systems.
In the case of exponential weight distribution, these queueing
systems correspond closely to totally asymmetric exclusion processes
or totally asymmetric zero-range processes. There are also 
close links with systems of non-colliding particles. 
See for example \cite{Neilsurvey} for a survey.

Of course, there are also strong connections between these 
models and random matrix theory. We mention one particular direction
related to the topic of this paper. 
Let $A_{n,k}$ be an $n\times k$ random matrix with i.i.d.\ entries,
and let $Y_{n,k} = A_{n,k}(A_{n,k})^*$. In the special
case of the Laguerre ensemble, where the 
common distribution of the entries 
is complex Gaussian, one has an explicit correspondence 
between the largest eigenvalue of $Y$ and 
the passage time to $(n,k)$ in a directed 
percolation model with exponential weights
(see for example Proposition 1.4 of \cite{Johshape},
and \cite{Doumerc} and Section 6.1 of \cite{BaiBenPec} for extensions).
For a general distribution, there may not exist 
an explicit mapping between the matrix
model and the directed percolation model, 
but we believe that a similar averaging 
mechanism to that observed in our context
will also play a role in the random matrix setting.
Thus on the basis of the analysis of the last passage time fluctuations
for paths close to the axis, we conjecture that
when $k$ and $n$ tend to infinity with suitable rates, the fluctuations
of the largest eigenvalue of $Y_{n,k}$ should depend only on the mean and the 
variance of the coefficients of $A_{n,k}$.

\section{Acknowledgments}
We thank Giambattista Giacomin, Yueyun Hu and Neil O'Connell for 
valuable discussions.

\end{document}